\documentclass[11pt, twoside, a4paper]{article}
\usepackage{makeidx}\usepackage{amssymb}

\newcommand{\prs}{\langle\;,\;\rangle}
\newcommand{\too}{\longrightarrow}
\newcommand{\om}{\omega}
\newcommand{\ch}{\check{e}}\newcommand{\tr}{{\mathrm{tr}}}

\newcommand{\ad}{{\mathrm{ad}}}

\newcommand{\G}{{\mathfrak{g}}}

\newcommand{\h}{{\cal H}}

\newcommand{\B}{{\cal B}}

\newcommand{\di}{\displaystyle}

\newcommand{\na}{\nabla}

\newcommand{\al}{\alpha}
\newcommand{\be}{\beta}
\newcommand{\ga}{\gamma}

\newcommand{\e}{\epsilon}

\newcommand{\la}{\lambda}

\newcommand{\de}{\delta}

\font\bb=msbm10

\def\B{\hbox{\bb B}}
\def\R{\hbox{\bb R}}

\def\N{\hbox{\bb N}}

\def\C{\hbox{\bb C}}

\newtheorem{theo}{Theorem}[section]
\newtheorem{pr}{Proposition}[section]
\newtheorem{Le}{Lemma}[section]
\newtheorem{co}{Corollary}[section]

\title{ Solutions of the  Yang-Baxter equations on orthogonal groups: the case of oscillator groups}

\author{Mohamed Boucetta-Alberto Medina}\date{}

\begin{document}\maketitle

\begin{abstract}
A Lie group is called orthogonal if it carries a bi-invariant pseudo-Riemannian metric. Oscillator Lie groups constitutes a  subclass of the class of orthogonal Lie groups.
 In this paper, we determine the Lie bialgebra structures and the solutions of the classical Yang-Baxter equation on a generic class of oscillator Lie groups. On the other hand, we show  that any solution of the classical Yang-Baxter equation on an orthogonal Lie group induces a flat left invariant pseudo-Riemannian metric in the dual Lie groups associated to this solution.
 This metric is geodesically complete if and only if the dual are unimodular. More generally, we show that any solution of the generalized Yang-Baxter equation on an orthogonal  Lie group determines a left invariant locally symmetric pseudo-Riemannian metric on the corresponding dual Lie groups. Applying this result to  oscillator Lie groups we get a large class of solvable Lie groups with flat left invariant Lorentzian metric.

\end{abstract}
\section{Introduction and main results}
The notion of Poisson-Lie group was first introduced by Drinfeld \cite{dr} and studied by Semenov-Tian-Shansky \cite{sts} (see also \cite{lu-wei}). It is known that every connected Poisson-Lie group arises from a Lie bialgebra and an important class of Lie bialgebras, the coboundary Lie bialgebras \cite{dr1}, are obtained by solving the generalized Yang-Baxter equation. Recall that a \emph{Lie bialgebra} is a Lie algebra $\G$ together with a linear map $\xi:\G\too\wedge^2\G$ such that:
\begin{enumerate}\item $\xi$ is a 1-cocycle with respect to the
adjoint action, i.e.,
\begin{equation}\label{cocycle}
\xi([u,v])=\ad_u\xi(v)-\ad_v\xi(u);\end{equation} \item the bracket
$[\;,\;]^*$ on the dual $\G^*$ given by
\begin{equation}\label{dualbracket}[\al,\be]^*(u)=\xi(u)(\al,\be),\quad u\in \G,
\al,\be\in\G^*\end{equation}satisfies the Jacobi identity.\end{enumerate}
A Lie bialgebra $(\G,\xi)$ is called \emph{coboundary Lie bialgebra} if there exists $r\in\wedge^2\G$ such that, for any $u\in\G$, $\xi(u)=\ad_{u}r$. In this case, the condition (\ref{cocycle}) is automatically satisfied and (\ref{dualbracket}) holds if and only if $r$ satisfies the \emph{generalized classical Yang-Baxter equation}:
\begin{equation}\label{gyb}\ad_u[r,r]=0,\quad\forall u\in\G,\end{equation}where $[r,r]\in\wedge^3\G$ is the Schouten bracket (see \cite{dz} for instance). A solution of the \emph{classical Yang-Baxter equation} is a bivector $r\in\wedge^2\G$ satisfying
\begin{equation}\label{cyb}[r,r]=0.\end{equation}
Bellavin and Drinfeld gave a  classification of the so-called quasi-triangular Lie bialgebras on simple complex Lie algebras (see \cite{Bel}). In \cite{delorme}, Delorme classified Manin triples (Lie bialgebras's equivalent) on reductive complex Lie algebras.
In \cite{z},  Poisson-Lie  structures on Heisenberg groups were classified.\\
In this paper, an oscillator group is a real simply connected Lie group which contains a Heisenberg group as a normal closed  subgroup of  codimension 1.
The four dimensional oscillator  group  has   its origin in the study of the  harmonic oscillator which is one of the most simple  non-relativist systems where the Schrodinger equation can be solved completely.    In  \cite{streater}, Streater described the representations of this group. Oscillator groups of dimension great than four
  have interesting features from the viewpoints of both Differential Geometry and Physics (see for instance \cite{gadea, levichev, levichev1, muller, nappi,nomura}).
 In \cite{medi-revoy} Medina proved that they are  the only non commutative simply connected solvable Lie groups which have a bi-invariant Lorentzian metric.

In what follows, we study  Poisson-Lie  structures  on oscillator groups. Indeed, we determine the Lie bialgebra structures and the  solutions of the generalized classical Yang-Baxter equation on a generic class of oscillator Lie algebras (see Theorems \ref{main1}-\ref{main2}).    We also show that any solution  of the classical Yang-Baxter equation on a Lie group endowed with a bi-invariant pseudo-Riemannian metric induces a flat left invariant pseudo-Riemannian metric on the dual Lie groups of $G$ associated to the solution (see Theorem \ref{flat}). Applying this result to oscillator Lie groups we get a large class of solvable Lie groups with flat left invariant Lorentzian metric. More generally, it is shown that any solution of the generalized Yang-Baxter equation on an orthogonal (or quadratic) Lie group determines a left invariant locally symmetric pseudo-Riemannian metric on the corresponding dual Lie groups.

  Note that our results imply interesting geometric properties on oscillator manifolds (quotient of oscillator groups by lattices) which constitute a large class since the conditions of existence of lattices on a given oscillator group can be fulfilled easily (see \cite{medi-revoy1}).

Let us now recall some facts and state our main results.

For $n\in\N^*$ and $\la=(\la_1,\ldots,\la_n)\in\R^n$ with
$0<\la_1\leq\ldots\leq\la_n$, the  $\la$-oscillator group, denoted by  $G_\la$, is the Lie  group which
the underlying manifold
$\R^{2n+2}=\R\times\R\times\C^n$ and product
$$(t,s,z).(t',s',z')=\left(t+t',s+s'+\frac12\sum_{j=1}^nIm\bar{z}_jexp(it\la_j)z_j',
\ldots,z_j+exp(it\la_j)z'_j,\ldots\right).$$

The Lie algebra of $G_\la$, denoted by  $\G_\la$, admits a basis
$\B=\left\{e_{-1},e_0,e_i,\ch_i,\right\}_{i=1,\ldots,n}$
where the brackets are given by
\begin{equation}\label{bracket}[e_{-1},e_j]=\la_j\ch_j,\qquad[e_{-1},\ch_j]
=-\la_je_j,\qquad[e_j,\ch_j]=e_0.\end{equation} The unspecified brackets are either zero
or given by antisymmetry.

We  call $G_\la$ or $\G_\la$ \emph{generic} if,  $0<\la_1<\ldots<\la_n$ and, for any $1\leq i<j<k\leq n$,  $\la_k\not=\la_i+\la_j$.

We shall denote by  $S$ the vector subspace of  $\G_\la$ spanned
by $\{e_i,\ch_i\}_{i=1,\ldots,n}$ and by  $\om$ the
 2-form on  $\G_\la$ given by
$$i_{e_{-1}}\om=i_{e_{0}}\om=0,\; \om(e_i,e_j)=\om(\ch_i,\ch_j)=0\quad\mbox{and}
\quad\om(e_i,\ch_j)=\de_{ij}.$$ The restriction of $\om$ to $S$ is a symplectic 2-form and, for any $u,v\in S$,
\begin{equation}\label{bracket1}[u,v]=\om(u,v)e_0.\end{equation}Moreover, $S$ is invariant by the derivation $\ad_{e_{-1}}$
and we have\begin{equation}\label{ad}\om(\ad_{e_{-1}}u,v)+\om(u,\ad_{e_{-1}}v)=0\quad
u,v\in S.\end{equation}
{\bf Notations.} For any $r_1,r_2\in\wedge^2\G_\la$, let $\om_{r_1,r_2}$ be the element of $\wedge^2\G_\la$ defined by
$$\om_{r_1,r_2}(\al,\be)=\frac12\left(\om(r_{1\#}(\al),r_{2\#}(\be))+\om(r_{2\#}(\al),r_{1\#}(\be)\right),$$where $r_{i\#}:\G_\la^*\too\G_\la$ is the endomorphism given by
$\be(r_{i\#}(\al))=r_i(\al,\be)$.  We denote (improperly) by $\wedge^2S$ the space of $r\in\wedge^2\G_\la$ satisfying $r_{\#}(e_{-1}^*)=r_{\#}(e_{0}^*)=0$ and   by $S^*$ the subspace of $\G_\la^*$ of  $\al$ such that $\al(e_{-1})=\al(e_0)=0$ and we define $e_{-1}^*,e_0^*\in\G_\la^*$  by ${e_{-1}^*}_{|S}={e_{0}^*}_{|S}=0$, $e_{-1}^*(e_0)=e_{0}^*(e_{-1})=0$ and $e_{-1}^*(e_{-1})=e_{0}^*(e_{0})=1$. We have clearly $\G_\la^*=\R e_{-1}^*\oplus\R e_{0}^*\oplus S^*.$ Finally,    for any endomorphism $J$ of $\G_\la$, we denote by $J^\dag$ the endomorphism of $\wedge^2\G_\la$, given by
$$J^\dag r(\al,\be)=r(J^*\al,\be)+r(\al,J^*\be),$$where $J^*:\G_\la^*\too\G_\la^*$ is the dual   of $J$.\\

We have:\begin{theo}\label{main1} Let $\G_\la$ be a generic oscillator Lie algebra. Then $\xi:\G_\la\too\wedge^2\G_\la$ defines a Lie bialgebra structure on $\G_\la$ if and only if there exists $r\in\wedge^2S$, $u_0\in S$ and a derivation $J:\G_\la\too\G_\la$ commuting with $\ad_{e_{-1}}$ and satisfying $J(e_{-1})=J(e_0)=0$ such that,
   for any $u\in\G_\la$, $$\xi(u)=\ad_u^\dag r+2e_0\wedge ((J+\ad_{u_0})(u)),$$ and
\begin{equation}\label{boucetta}\om_{r,\ad_{e_{-1}}^\dag r}-(J^\dag\circ \ad_{e_{-1}}^\dag)  r=0.\end{equation}
Moreover, in this case, the Lie bracket on $\G_\la^*$ defined by (\ref{dualbracket}) is given by
\begin{equation}\label{bracketmain}\left\{\begin{array}{l}
\;[e_{0}^*,\al]^*=2J^*\al-2(\ad^*_{e_{-1}}\al)(u_0)e_{-1}^*+ i_{r_{\#}(\al)}\om,\\
\;[\al,\be]^*=\ad_{e_{-1}}^\dag r(\al,\be)e_{-1}^*,\end{array}\right.\end{equation}where  $\al,\be\in S^*$ and $e_{-1}^*$ is a central element.
 \end{theo}
From the expression of the brackets above, we can deduce immediately the following result.
\begin{co}\label{co}Let $\G_\la$ be a generic oscillator Lie algebra and $\xi(u)=
\ad_u^\dag r+2e_0\wedge ((J+\ad_{u_0})(u))$  a Lie bialgebra structure on $\G_\la$. Denote by $2p$ the dimension of the kernel of the restriction of $\ad_{e_{-1}}^\dag r$ to $S^*$. Then $(\G^*_\la,[\;,\;]^*)$ is isomorphic to the semi-direct product of $\R e_{0}^*$ by the ideal $\h_{2(n-p)+1}\oplus\R^{2p}$ where $\h_{2(n-p)+1}$ is the $2(n-p)+1$-dimensional Heisenberg Lie algebra, $\R^{2p}$ the $2p$-dimensional Abelian Lie algebra and the action of $e_{0}^*$ is given by the first relation in (\ref{bracketmain}). In particular, $(\G^*_\la,[\;,\;]^*)$ is solvable. Moreover, $(\G^*_\la,[\;,\;]^*)$
is unimodular iff $\sum_{i=1}^nr(e_i,\ch_i)=0$.\end{co}

\begin{theo}\label{main2}Let $\G_\la$ be a generic oscillator Lie algebra. Then:
\begin{enumerate}\item A bivector $r\in\wedge^2\G_\la$ is a solution of the generalized  Yang-Baxter equation if and only if there exists $u_0\in S$, $r_{0}\in\wedge^2S$ and $\al\in\R$ such that $r=2\al e_0\wedge e_{-1}+e_0\wedge u_0+r_0$
and\begin{equation}\label{gyb1}\om_{r_0,\ad_{e_{-1}}^\dag r_0}+\al (\ad_{e_{-1}}^\dag\circ \ad_{e_{-1}}^\dag)r_0=0.\end{equation}
\item A bivector $r\in\wedge^2\G_\la$ is a solution of the  classical Yang-Baxter equation if and only if there exists $u_0\in S$, $r_{0}\in\wedge^2S$ and $\al\in\R$ such that $r=\al e_0\wedge e_{-1}+e_0\wedge u_0+r_0$
and\begin{equation}\label{cyb1}\om_{r_0,r_0}+\al \ad_{e_{-1}}^\dag r_0=0.\end{equation}
\end{enumerate}
\end{theo}
There are some comments on Theorems \ref{main1}-\ref{main2}:
\begin{enumerate}\item  Theorems \ref{main1}-\ref{main2} reduce the problem of finding Lie bialgebras structures or solutions of Yang-Baxter equations on a generic oscillator Lie algebra to solving  (\ref{boucetta}), (\ref{gyb1}) and  (\ref{cyb1}). Or these equations involve only the symplectic space $(S,\om)$ and the restrictions of the derivations $J$ and $\ad_{e_{-1}}$ to  $S$. Note that the space of derivations of $\G_\la$ satisfying $J(e_0)=J(e_{-1})=0$ is isomorphic to the space of endomorphisms of $S$ skew-symmetric with respect to $\om$.

\item In Section \ref{example}, we give the solutions of (\ref{boucetta}), (\ref{gyb1}) and  (\ref{cyb1}) when $\dim\G_\la\leq 6$ (see Propositions \ref{dim4}-\ref{dim6}). To solve those equations in the general case is very difficult, however we will give in Section \ref{example} a large class of solutions.

    \item On a general oscillator algebra (not necessary generic), the 1-cocyles satisfying the conditions of Theorem \ref{main1} (resp. the $r$ satisfying the conditions of Theorem \ref{main2}) defines a large class of Lie bialgebras structures on $\G_\la$ (resp.  a large class of solutions of the generalized Yang-Baxter equation on $\G_\la$).
        \item When $\al=0$, $r_0$ is a solution of (\ref{cyb1}) if and only if $\mbox{Im}r_{0\#}$ is an
        even dimensional $\om$-isotropic subspace of $S$.\\ Conversely, let $F$ be an even dimensional $\om$-isotropic subspace of $S$. Choose a nondegenerate skew-symmetric 2-form $\mu$ on $F$ and define $r_0\in\wedge^2\G_\la$ by
        $$r_0(\al,\be)=\be\left(\mu_{\#}\circ i^*(\al)\right),$$where $i:F\too\G_\la$ is canonical inclusion and $\mu_{\#}:F^*\too F$ the isomorphism associated to $\mu$. One can check easily that $r_0$ is a solution of (\ref{cyb1}).

\end{enumerate}

The second part of our study involves bi-invariant pseudo-Riemannian metrics on Lie groups and solutions of the generalized  Yang-Baxter equation.

Let $(G,k)$ be a Lie group endowed with a bi-invariant pseudo-Riemannian metric. The value of $k$ at identity induces on the Lie algebra $\G$ of $G$ an adjoint invariant non degenerate bilinear symmetric form $\prs$. Such a Lie algebra
 is called an \emph{orthogonal} (or \emph{quadratic}) Lie algebra. For instance, semi-simple Lie algebras and oscillator Lie algebras are orthogonal (see
\cite{medi-revoy}). \\ Let $r$ be a solution of (\ref{gyb}) on an orthogonal Lie algebra $(\G,\prs)$. Then $r$ defines on $\G^*$ a Lie bracket by
 \begin{equation}\label{dualbracket1}[\al,\be]_r=\ad_{r_{\#}(\be)}^*\al-\ad_{r_{\#}(\al)}^*\be.\end{equation}Consider the bilinear form on $\G^*$ given by
 $$\langle\al,\be\rangle^*=\langle \phi^{-1}(\al),\phi^{-1}(\be)\rangle,$$where $\phi(u)=\langle u,.\rangle.$\\
 Let us denote by $G^*_r$ a Lie group with Lie algebra $(\G^*,[\;,\;]_r)$, by $k^*$ the left invariant pseudo-Riemannian metric whose value at the identity is $\prs^*$ and by $\na^*$ its Levi-Civita connexion. With the notations above, we have the following result.
\begin{theo}\label{flat} Let $(G,k)$ be a Lie group endowed with a bi-invariant pseudo-Riemannian metric and let
$r$ be a solution of (\ref{gyb})  on $G$. Then:
\begin{enumerate}
\item  $(G^*_r,k^*)$ is a locally symmetric pseudo-Riemannian manifold, i.e., $$\na^*R=0,$$ where $R$ is the curvature of $k^*$. In particular, $R$ vanishes identically when $r$ is a solution of (\ref{cyb}).

\item  If $k^*$ is flat then it  is complete if and only if $G_r^*$ is unimodular and in this case $G_r^*$ is solvable.

    \end{enumerate}\end{theo}

  In \cite{mr} Theorem 3.9, Medina an Revoy obtained a similar result as Theorem \ref{flat} when $r$ is an invertible solution of (\ref{cyb}). For more details on left invariant metrics on quadratic Lie groups see \cite{brmedi1}. On the other hand, the results above complete the results obtained by Bordemann in \cite{bordemann}.

There are some consequences of Theorem \ref{flat}. As above $(G,k)$ is a Lie group endowed with a bi-invariant pseudo-Riemannian metric, $(\G,\prs)$ its associated orthogonal Lie algebra and $r\in\wedge^2\G$.
\begin{enumerate}

\item If $r$ is a solution of (\ref{cyb})  and $k$ is definite (for instance, $k$ is the Killing form of a semi-simple compact Lie group) then, according to a result of Milnor \cite{milnor}, $\G^*$ is a semi-direct product of two Abelian  Euclidian Lie algebras,  one of them acting on the other by infinitesimal isometries.

    \item If $r$ is a solution of (\ref{cyb}), the sequence $$0\too\ker r_{\#}\too\G^*\too\mbox{Im}r_{\#}\too0$$ is an exact sequence of Lie algebras with $\ker r_{\#}$ is an abelian ideal of $\G^*$ and $\mbox{Im}r_{\#}$ is a symplectic Lie algebra. In particular if $G$ is solvable  then $G_r^*$ is solvable. The sequence is also a natural exact sequence of left symmetric algebras (see \cite{mr}).

      \item   We can paraphrase a part of Theorem \ref{flat} as follows: Any symplectic Lie subgroup $S$ of an ortogonal Lie group  $(G,k)$ determines an unique simply connected Lie group
$G*$   such that the left invariant pseudo-Riemannian metric $k*$ is flat.
For example, the affine group  $S=\rm{Aff}( n ,\R )$ is a symplectic Lie subgroup of the orthogonal  Lie group $G= \rm{GL} ( n+1,\R )$ . Hence there is a simply connected Lie group dual of $G$ endowed with a flat left invariant peudo-Riemannian metric non complete.

        \item Any solution of (\ref{gyb}) (resp. (\ref{cyb})) on a reductive or semi-simple Lie group gives  rise to a locally symmetric (resp. flat) left invariant pseudo-Riemannian metric on the dual groups.
            \item Recall that a Lie group can carry many several orthogonal structures non isomorphic. This is the case, for example, if the groupe is  semi-simple an non simple. Consequently the dual Lie groups of $G$  can be have many
flat left invariant pseudo-Riemannian metrics.

\item  As shown by Medina in \cite{medi-revoy}, oscillator Lie algebras are orthogonal. Indeed, for $x\in\G_\la$, let
 $$x=x_{-1}e_{-1}+x_0e_0+\sum_{i=1}^n\left(x_ie_i+\check{x}_i\check{e}_i\right).$$
 The nondegenerate quadratic form
 \begin{equation}\label{kl}\textbf{k}_\la(x,x):=2x_{-1}x_0+\sum_{i=1}^n\frac1{\la_i}(x_i^2+\check{x}_i^2)\end{equation} defines  a
 Lorentzian bi-invariant metric on $G_\la$. According to Theorem \ref{flat}, any solution of (\ref{gyb}) (resp. (\ref{cyb})) on $G_\la$  gives  rise to a locally symmetric (resp. flat) left invariant Lorentzian metric on the dual group.

 \end{enumerate}

 The paper is organized as follows: Section \ref{bialgebra} is devoted to the study of Lie bialgebras structures on oscillator Lie algebras and culminates by proving  Theorems \ref{main1}-\ref{flat}. In Section \ref{example}, we give the solutions of (\ref{boucetta}), (\ref{gyb1}) and  (\ref{cyb1}) when $\dim\G_\la\leq 6$. In section \ref{example1}, we illustrate Theorem  \ref{flat} by building   an example of a complete 6-dimensional Lie group endowed with a left invariant flat Lorentzian metric.  Moreover, since the Killing form of $SL(2,\R)$ is non degenerate Lorentzian,  we give the solutions of (\ref{cyb}) on $SL(2,\R)$ and we give the associated flat Lorentzian dual groups.

\section{Lie bialgebras structures on oscillator Lie algebras}\label{bialgebra}
 This section is devoted to the study of bialgebras structures on oscillator Lie algebras. It culminates with a proof of Theorems \ref{main1}-\ref{flat}.\\ Let
$\B^*=\left\{e_{-1}^*,e_0^*,e_i^*,\ch_i^*\right\}_{i=1,\ldots,n}$
be the dual basis of
$\B$.
From (\ref{bracket}), we get that the non vanishing $\ad_u^*\al$ with $u\in\B$ and $\al\in\B^*$ are
\begin{eqnarray}\label{adjoint}
\ad_{e_{-1}}^*e_i^*&=&-\la_i\ch_i^*,\;
\ad_{e_{-1}}^*\ch_i^*=\la_ie_i^*,\nonumber\\
\ad_{e_{i}}^*e_0^*&=&\ch_i^*,\;\ad_{e_{i}}^*\ch_i^*=-\la_ie_{-1}^*,\label{adjoint}\\
\ad_{\ch_{i}}^*e_0^*&=&-e_i^*,\;\ad_{\ch_{i}}^*e_i^*=\la_ie_{-1}^*.\nonumber\end{eqnarray}
An important class of 1-cocycles are the coboundaries, i.e, $\xi:\G_\la\too\G_\la\wedge\G_\la$ such that, for any $u\in\G_\la$, $\xi(u)=\ad_u^\dag r$ where $r\in\G_\la\wedge\G_\la$. In this case, the bracket defined by (\ref{dualbracket}) is given by
$$[\al,\be]^*=\ad_{r_{\#}(\be)}^*\al-\ad_{r_{\#}(\al)}^*\be,\quad\al,\be\in\G_\la^*.$$ In the proof of Theorem \ref{main1} we need to compute the expression of this bracket in the basis $\B^*$. A direct computation using (\ref{adjoint}) gives:
\begin{equation}\label{exactbracket}\left\{\begin{array}{lll}
\;[e_{-1}^*,\al]^*&=&-\ad_{r_{\#}(e_{-1}^*)}^*\al,\quad\al\in\G_\la^*,\\
\;[e_0^*,e_i^*]^*&=&-\la_ir(e_0^*,\ch_i^*)e_{-1}^*+\la_ir(e_0^*,e_{-1}^*)\ch_i^*+
\di\sum_{j=1}^n\left(r(e_i^*,e_j^*)\ch_j^*-r(e_i^*,\ch_j^*)e_j^*\right),\\
\;[e_0^*,\ch_i^*]^*&=&\la_ir(e_0^*,e_i^*)e_{-1}^*-\la_ir(e_0^*,e_{-1}^*)e_i^*+\di\sum_{j=1}^n\left(r(\ch_i^*,e_j^*)\ch_j^*-
r(\ch_i^*,\ch_j^*)e_j^*\right),\\
\;[e_i^*,e_j^*]^*&=&\left(-\la_ir(\ch_i^*,e_j^*)-\la_jr(e_i^*,\ch_j^*)\right)e_{-1}^*+\la_ir(e_{-1}^*,e_j^*)\ch_i^*+
\la_jr(e_{i}^*,e_{-1}^*)\ch_j^*,\\
\;[e_i^*,\ch_j^*]^*&=&\left(-\la_ir(\ch_i^*,\ch_j^*)+\la_jr(e_i^*,e_j^*)\right)e_{-1}^*+\la_ir(e_{-1}^*,\ch_j^*)\ch_i^*-
\la_jr(e_{i}^*,e_{-1}^*)e_j^*,\\
\;[\ch_i^*,\ch_j^*]^*&=&\left(\la_ir(e_i^*,\ch_j^*)+\la_jr(\ch_i^*,e_j^*)\right)e_{-1}^*-\la_ir(e_{-1}^*,\ch_j^*)e_i^*-
\la_jr(\ch_{i}^*,e_{-1}^*)e_j^*.\end{array}\right.\end{equation}

{\bf Notation.} For any linear map $\xi:\G_\la\too\G_\la\wedge\G_\la$ and for any  $i,j=1,\ldots,n$, we denote by $\al_{-1,0},\; \al_{-1,i},\; \check{\al}_{-1,i},\; a_{0,i},\; \check{a}_{0,i},\; b_{i,j},\; \check{b}_{i,j},\; c_{i,j}$ the elements of $\G_\la^*$ whose values at $u\in\G_\la$ are given by
\begin{eqnarray*}
\al_{-1,0}(u)&=&\xi(u)(e_{-1}^*,e_0^*),\;\al_{-1,i}(u)=\xi(u)(e_{-1}^*,e_i^*),\;
\check{\al}_{-1,i}(u)=\xi(u)(e_{-1}^*,\ch_i^*),\\
a_{0,i}(u)&=&\xi(u)(e_{0}^*,e_i^*),\;
\check{a}_{0,i}(u)=\xi(u)(e_{0}^*,\ch_i^*),\\
b_{i,j}(u)&=&\xi(u)(e_{i}^*,e_j^*),\;
\check{b}_{i,j}(u)=\xi(u)(\ch_{i}^*,\ch_j^*),\;
c_{i,j}(u)=\xi(u)(e_{i}^*,\ch_j^*).\end{eqnarray*}

\begin{pr}\label{pr1} Let $\xi:\G_\la\too\G_\la\wedge\G_\la$ be a 1-cocycle
with respect to the adjoint action. Then $\xi(e_0)=0$.\end{pr}

{\bf Proof.} Since $e_0$ is a central element in $\G_\la$, the
cocycle condition implies
$\ad_u\xi(e_0)=0,$ for any $u\in\G_\la$.  Hence, by using (\ref{adjoint}), we get, for $i,j=1,\ldots,n$,
 \begin{eqnarray*}
 0&=&\ad_{e_i}\xi(e_0)(e_0^*,\ch_i^*)=\xi(e_0)(\ad^*_{e_i}e_0^*,\ch_i^*)+
 \xi(e_0)(e_0^*,\ad^*_{e_i}\ch_i^*)=-\la_i\xi(e_0)(e_0^*,e_{-1}^*),\\
0&=&\ad_{e_i}\xi(e_0)(e_0^*,e_{-1}^*)=\xi(e_0)(\ad_{e_i}^*e_0^*,e_{-1}^*)=
\xi(e_0)(\ch_i^*,e_{-1}^*),\\
0&=&\ad_{\ch_i}\xi(e_0)(e_0^*,e_{-1}^*)=\xi(e_0)(\ad_{\ch_i}^*e_0^*,e_{-1}^*)=
-\xi(e_0)(e_i^*,e_{-1}^*),\\
0&=&\ad_{e_i}\xi(e_0)(e_0^*,e_j^*)=\xi(e_0)(\ad_{e_i}^*e_0^*,e_j^*)=
\xi(e_0)(\ch_i^*,e_j^*),\\
0&=&\ad_{e_{-1}}\xi(e_0)(e_0^*,e_i^*)=\xi(e_0)(e_0^*,\ad_{e_{-1}}^*e_i^*)=
-\la_i\xi(e_0)(e_0^*,\ch_i^*),\\
0&=&\ad_{e_{-1}}\xi(e_0)(e_0^*,\ch_i^*)=\xi(e_0)(e_0^*,\ad_{e_{-1}}^*\ch_i^*)=
\la_i\xi(e_0)(e_0^*,e_i^*),\\
0&=&\ad_{e_i}\xi(e_0)(e_0^*,\ch_j^*)=\xi(e_0)(\ad_{e_i}^*e_0^*,\ch_j^*)+
\xi(e_0)(e_0^*,\ad_{e_i}^*\ch_j^*)=
\xi(e_0)(\ch_i^*,\ch_j^*),\\
0&=&\ad_{\ch_i}\xi(e_0)(e_0^*,e_j^*)=\xi(e_0)(\ad_{\ch_i}^*e_0^*,e_j^*)+
\xi(e_0)(e_0^*,\ad_{\ch_i}^*e_j^*)=- \xi(e_0)(e_i^*,e_j^*),
\end{eqnarray*}which shows that $\xi(e_0)=0$.\hfill $\square$\\

\begin{pr}\label{pr2}Let $\xi:\G_\la\too\G_\la\wedge\G_\la$ be a 1-cocycle
with respect to the adjoint action. Then, for any $i,j=1,\ldots,n$ with $i\not=j$,
\begin{eqnarray*}
b_{i,j}(e_i)&=&b_{i,j}(e_j)=0,\\
c_{j,j}(\ch_i)&=&c_{j,j}(e_i)=c_{i,j}(e_i)=c_{i,j}(\ch_j)=c_{i,i}(e_{-1})=0,\\
\check{b}_{i,j}(\ch_i)&=&\check{b}_{i,j}(\ch_j)=0.\end{eqnarray*}

\end{pr}

{\bf Proof.} Let $i,j=1,\ldots,n$ with $i\not=j$. The cocycle condition, Proposition \ref{pr1}, (\ref{bracket}) and (\ref{adjoint}) imply
\begin{eqnarray*}
0&=&\ad_{\ch_i}\xi(\ch_j)(e_0^*,\ch_j^*)-\ad_{\ch_j}\xi(\ch_i)(e_0^*,\ch_j^*)=
-\xi(\ch_j)(e_i^*,\ch_j^*)+\xi(\ch_i)(e_j^*,\ch_j^*)\\
&=&-c_{i,j}(\ch_j)+c_{j,j}(\ch_i).\\
0&=&\ad_{\ch_i}\xi(e_j)(e_0^*,e_j^*)-\ad_{e_j}\xi(\ch_i)(e_0^*,e_j^*)=
-\xi(e_j)(e_i^*,e_j^*)-\xi(\ch_i)(\ch_j^*,e_j^*)\\
&=&-b_{i,j}(e_j)+c_{j,j}(\ch_i).\\
0&=&\ad_{e_j}\xi(\ch_j)(e_0^*,e_i^*)-\ad_{\ch_j}\xi(e_j)(e_0^*,e_i^*)=
\xi(\ch_j)(\ch_j^*,e_i^*)+\xi(e_j)(e_j^*,e_i^*)\\
&=&-c_{i,j}(\ch_j)-b_{i,j}(e_j).
\end{eqnarray*}We deduce that
$b_{i,j}(e_i)=c_{j,j}(\ch_i)=c_{i,j}(\ch_j)=0.$
On the other hand,
\begin{eqnarray*}
0&=&\ad_{e_i}\xi(e_j)(e_0^*,e_j^*)-\ad_{e_j}\xi(e_i)(e_0^*,e_j^*)=
\xi(e_j)(\ch_i^*,e_j^*)-\xi(e_i)(\ch_j^*,e_j^*)\\
&=&-c_{j,i}(e_j)+c_{j,j}(e_i).\\
0&=&\ad_{e_i}\xi(\ch_j)(e_0^*,\ch_j^*)-\ad_{\ch_j}\xi(e_i)(e_0^*,\ch_j^*)=
\xi(\ch_j)(\ch_i^*,\ch_j^*)+\xi(e_i)(e_j^*,\ch_j^*)\\
&=&\check{b}_{i,j}(\ch_j)+c_{j,j}(e_i).\\
0&=&\ad_{\ch_j}\xi(e_j)(e_0^*,\ch_i^*)-\ad_{e_j}\xi(\ch_j)(e_0^*,\ch_i^*)=
-\xi(e_j)(e_j^*,\ch_i^*)-\xi(\ch_j)(\ch_j^*,\ch_i^*)\\
&=&-c_{j,i}(e_j)+\check{b}_{i,j}(\ch_j).
\end{eqnarray*}
Thus
$\check{b}_{i,j}(\ch_i)=c_{j,j}(e_i)=c_{i,j}(e_i)=0.$  To complete the proof, we need to show that
$c_{i,i}(e_{-1})=0$. Indeed, by applying (\ref{cocycle}), respectively, to $(e_{-1},e_i)$ and $(e_{-1},\ch_i)$, we get\begin{eqnarray*}
\la_i\xi(\ch_i)(e_0^*,e_i^*)&=&\ad_{e_{-1}}\xi(e_i)(e_0^*,e_i^*)
                               -\ad_{e_i}\xi(e_{-1})(e_0^*,e_i^*)\\
 &=&-\la_i\xi(e_i)(e_0^*,\ch_i^*)-\xi(e_{-1})(\ch_i^*,e_i^*),\\
 -\la_i\xi(e_i)(e_0^*,\ch_i^*)&=&\ad_{e_{-1}}\xi(\ch_i)(e_0^*,\ch_i^*)
                               -\ad_{\ch_i}\xi(e_{-1})(e_0^*,\ch_i^*)\\
 &=&\la_i\xi(\ch_i)(e_0^*,e_i^*)+\xi(e_{-1})(e_i^*,\ch_i^*).\end{eqnarray*}
 Hence
 $$c_{i,i}(e_{-1})=\la_i(a_{0,i}(\ch_i)+\check{a}_{0,i}(e_i))=-\la_i(\check{a}_{0,i}(e_i)+a_{0,i}(\ch_i)),$$
and then $c_{i,i}(e_{-1})=0$.
\hfill $\square$

\begin{pr}\label{pr3}Let $\xi:\G_\la\too\G_\la\wedge\G_\la$ be a 1-cocycle
on $\G_\la$ with respect to the adjoint action and let $k,l=1,\ldots,n$. Then, for any $1\leq j\leq n$ such that $j\not=k$ and $j\not=l$, we have:
\begin{enumerate} \item
if $\la_j\not=\la_k+\la_l$ and $\la_j\not=|\la_l-\la_k|$, then
$$b_{k,l}(e_j)=b_{k,l}(\ch_j)=\check{b}_{k,l}(e_j)=\check{b}_{k,l}(\ch_j)=c_{k,l}(e_j)=c_{k,l}(\ch_j)=0;$$
 \item if
$\la_j=\la_k+\la_l$, then there exists $a,b\in\R$ such that
\begin{eqnarray*}
(b_{k,l}(e_j),c_{k,l}(e_j),c_{l,k}(e_j),
\check{b}_{k,l}(e_j))&=&(-a,b,-b,a)\\
(b_{k,l}(\ch_j),c_{k,l}(\ch_j),c_{l,k}(\ch_j),
\check{b}_{k,l}(\ch_j))&=&(-b,-a,a,b);\end{eqnarray*}

 \item If
$\la_j=|\la_l-\la_k|$, then there exists $a,b\in\R$ such that
\begin{eqnarray*}
(b_{k,l}(e_j),c_{k,l}(e_j),c_{l,k}(e_j),
\check{b}_{k,l}(e_j))&=&(a,-b,-b,a)\\
(b_{k,l}(\ch_j),c_{k,l}(\ch_j),c_{l,k}(\ch_j),
\check{b}_{k,l}(\ch_j))&=&(b,\e a,\e a,b),\end{eqnarray*}where $\e$ is the signe of $\la_l-\la_k$.
\end{enumerate}
\end{pr}

{\bf Proof.} By applying (\ref{cocycle}) to $(e_{-1},e_j)$ and $(e_{-1},\ch_j)$, we get
\begin{eqnarray*}\ad_{
e_{-1}}\xi(e_j )-\ad_{ e_j}\xi(e_{-1}
)&=&\la_j\xi(\ch_j)\\
\ad_{ e_{-1}}\xi(\ch_j )-\ad_{ \ch_j}\xi(e_{-1}
)&=&-\la_j\xi(e_j).\end{eqnarray*} By evaluating these relations, respectively, on
$(e_k^*,e_l^*),(e_k^*,\ch_l^*),(\ch_k^*,e_l^*),(\ch_k^*,\ch_l^*),$ we deduce
\begin{eqnarray*}
\la_j\xi(\ch_j)(e_k^*,e_l^*)&=&-\la_k\xi(e_j)(\ch_k^*,e_l^*)-\la_l\xi(e_j)(e_k^*,\ch_l^*),\\
\la_j\xi(\ch_j)(e_k^*,\ch_l^*)&=&-\la_k\xi(e_j)(\ch_k^*,\ch_l^*)+\la_l\xi(e_j)(e_k^*,e_l^*),\\
\la_j\xi(\ch_j)(\ch_k^*,e_l^*)&=&\la_k\xi(e_j)(e_k^*,e_l^*)-\la_l\xi(e_j)(\ch_k^*,\ch_l^*),\\
\la_j\xi(\ch_j)(\ch_k^*,\ch_l^*)&=&\la_k\xi(e_j)(e_k^*,\ch_l^*)+\la_l\xi(e_j)(\ch_k^*,e_l^*),\\
-\la_j\xi(e_j)(e_k^*,e_l^*)&=&-\la_k\xi(\ch_j)(\ch_k^*,e_l^*)-
\la_l\xi(\ch_j)(e_k^*,\ch_l^*),\\
-\la_j\xi(e_j)(e_k^*,\ch_l^*)&=&-\la_k\xi(\ch_j)(\ch_k^*,\ch_l^*)+
\la_l\xi(\ch_j)(e_k^*,e_l^*),\\
-\la_j\xi(e_j)(\ch_k^*,e_l^*)&=&\la_k\xi(\ch_j)(e_k^*,e_l^*)-
\la_l\xi(\ch_j)(\ch_k^*,\ch_l^*),\\
-\la_j\xi(e_j)(\ch_k^*,\ch_l^*)&=&\la_k\xi(\ch_j)(e_k^*,\ch_l^*)
+\la_l\xi(\ch_j)(\ch_k^*,e_l^*).\end{eqnarray*}

These relations can by written
\begin{equation}\label{matrix1}B{X}=-\la_j\check{X}\quad\mbox{and}\quad B\check{X}=\la_jX.\end{equation}where
$$B=\left(\begin{array}{cccc}0&\la_l&\la_k&0\\-\la_l&0&0&\la_k\\-\la_k&0&0&\la_l\\
0&-\la_k&-\la_l&0\end{array}\right),$$
$
X=(b_{k,l}(e_j),c_{k,l}(e_j),-c_{l,k}(e_j),
\check{b}_{k,l}(e_j))$ and $
\check{X}=(b_{k,l}(\ch_j),c_{k,l}(\ch_j),-c_{l,k}(\ch_j),
\check{b}_{k,l}(\ch_j)).$

The equation (\ref{matrix1}) implies that
$$B^2X=-\la_j^2 X\quad\mbox{and}\quad
B^2\check{X}=-\la_j^2\check{X}.$$ Now the eigenvalues of $B^2$ are
$-(\la_l-\la_k)^2$ and $-(\la_l+\la_k)^2$ and the corresponding
eigenspaces are $span\{(1,0,0,1),(0,-1,1,0)\}$ and
$span\{(-1,0,0,1),(0,1,1,0)\}$. Hence
\begin{enumerate}\item If $\la_j^2\not=(\la_l+\la_k)^2$ and $\la_j^2\not=(\la_l-\la_k)^2$
then $X=\check{X}=0$. \item If $\la_j=|\la_l-\la_k|$ then
\begin{eqnarray*}
X&=&a(1,0,0,1)+b(0,-1,1,0),\\
\check{X}&=&\check{a}(1,0,0,1)+\check{b}(0,-1,1,0)\end{eqnarray*}and
(\ref{matrix1}) is equivalent to
$(\la_l-\la_k)a=-\la_j\check{b}\quad\mbox{and}\quad
(\la_k-\la_l)b=-\la_j\check{a}$ and we get the desired relation.

 \item If $\la_j=\la_l+\la_k$ then
\begin{eqnarray*}
X&=&a(-1,0,0,1)+b(0,1,1,0),\\
\check{X}&=&\check{a}(-1,0,0,1)+\check{b}(0,1,1,0)\end{eqnarray*}and
(\ref{matrix1}) is equivalent to
$a=-\check{b}\quad\mbox{and}\quad
b=\check{a}$ and we get the desired relation. \hfill $\square$
\end{enumerate}

\begin{pr}\label{pr4} Let $\xi:\G_\la\too\G_\la\wedge\G_\la$ be a 1-cocycle
on $\G_\la$ with respect to the adjoint action. Then
 for any $1\leq i,j\leq n$,
 $$\left\{\begin{array}{ccl}
 \al_{-1,0}(e_{-1})&=&0,\\
 c_{i,i}(e_{-1})&=&0,\\
 a_{0,i}(e_i)&=&\check{a}_{0,i}(\ch_i),\\
a_{0,i}(\ch_i)&=&-\check{a}_{0,i}(e_i),\\
b_{i,j}(e_{-1})&=&\la_j\check{a}_{0,j}(\ch_i)-\la_ia_{0,j}(e_i),\\
\check{b}_{i,j}(e_{-1})&=&\la_ja_{0,j}(e_i)-\la_i\check{a}_{0,j}(\ch_i),\\
c_{j,i}(e_{-1})&=&\la_ia_{0,j}(\ch_i)+\la_j\check{a}_{0,j}(e_i),\\
c_{i,j}(e_{-1})&=&-\la_i\check{a}_{0,j}(e_i)-\la_ja_{0,j}(\ch_i).\end{array}\right.$$

\end{pr}

 {\bf Proof.} Let $i\not=j$. By applying (\ref{cocycle}) to $(e_{-1},e_i)$ and $(e_{-1},\ch_i)$, we get\begin{eqnarray*}\ad_{
e_{-1}}\xi(e_i )-\ad_{ e_i}\xi(e_{-1}
)&=&\la_i\xi(\ch_i),\\
\ad_{ e_{-1}}\xi(\ch_i )-\ad_{ \ch_i}\xi(e_{-1}
)&=&-\la_i\xi(e_i).\end{eqnarray*} The evaluation of these two relations,  respectively, on
$(e_0^*,e_i^*),(e_0^*,e_j^*),(e_0^*,\ch_i^*),(e_0^*,\ch_j^*)$ gives
$$\left\{\begin{array}{ccl}
c_{i,i}(e_{-1})&=&\la_i(a_{0,i}(\ch_i)+\check{a}_{0,i}(e_i)),\\
c_{j,i}(e_{-1})&=&\la_ia_{0,j}(\ch_i)+\la_j\check{a}_{0,j}(e_i),\\
\al_{-1,0}(e_{-1})&=&a_{0,i}(e_i)-\check{a}_{0,i}(\ch_i) ,\\
\check{b}_{i,j}(e_{-1})&=&\la_ja_{0,j}(e_i)-\la_i\check{a}_{0,j}(\ch_i),\\
\al_{-1,0}(e_{-1})&=&\check{a}_{0,i}(\ch_i)-a_{0,i}(e_i),\\
b_{i,j}(e_{-1})&=&\la_j\check{a}_{0,j}(\ch_i)-\la_ia_{0,j}(e_i),\\
c_{i,i}(e_{-1})&=&-\la_i(\check{a}_{0,i}(e_i)+a_{0,i}(\ch_i)),\\
c_{i,j}(e_{-1})&=&-\la_i\check{a}_{0,j}(e_i)-\la_ja_{0,j}(\ch_i).\end{array}\right.$$
and the proposition follows.\hfill $\square$

\begin{Le}\label{centrallemma} Let $\xi:\G_\la\too\G_\la\wedge\G_\la$ be a bialgebra
structure on $\G_\la$. Then $e_{-1}^*$ is a central element of the dual Lie algebra $\G_\la^*$.\end{Le}

{\bf Proof.} Remark that $e_{-1}^*$ is
a central element if and only if $\al_{-1,0}=\al_{-1,i}=\check{\al}_{-1,i}=0$ for $i=1,\ldots,n$.

Note first that, according to Proposition \ref{pr1}, $\al_{-1,0}(e_0)=\al_{-1,i}(e_0)=\check{\al}_{-1,i}(e_0)=0$ for $i=1,\ldots,n$.

Let $1\leq i,j,k\leq n$ such that $i\not=j$. By using the 1-cocycle condition, (\ref{bracket}) and (\ref{adjoint}) we get
\begin{eqnarray*}
0&=&\ad_{e_i}\xi(e_j)(\ch_i^*,e_k^*)-\ad_{e_j}\xi(e_i)(\ch_i^*,e_k^*)=
-\la_i\xi(e_j)(e_{-1}^*,e_k^*)=-\la_i\al_{-1,k}(e_j),\\
0&=&\ad_{e_i}\xi(\ch_j)(\ch_i^*,e_i^*)-\ad_{\ch_j}\xi(e_i)(\ch_i^*,e_i^*)=
-\la_i\xi(\ch_j)(e_{-1}^*,e_i^*)=-\la_i\al_{-1,i}(\ch_j),\\
0&=&\ad_{\ch_i}\xi(\ch_j)(\ch_k^*,e_i^*)-\ad_{\ch_j}\xi(\ch_i)(\ch_k^*,e_i^*)=
\la_i\xi(\ch_j)(\ch_k^*,e_{-1}^*)=-\la_i\check{\al}_{-1,k}(\ch_j),\\
0&=&\ad_{\ch_i}\xi(e_j)(\ch_i^*,e_i^*)-\ad_{e_j}\xi(\ch_i)(\ch_i^*,e_i^*)=
\la_i\xi(e_j)(\ch_i^*,e_{-1}^*)=-\la_i\check{\al}_{-1,i}(e_j).
\end{eqnarray*}
Thus $\al_{-1,k}(e_j)=\al_{-1,i}(\ch_j)=\check{\al}_{-1,k}(\ch_j)=\check{\al}_{-1,i}(e_j)=0.$

On the other hand,
\begin{eqnarray*}
0&=&\ad_{e_j}\xi(e_i)(\ch_j^*,\ch_i^*)-\ad_{e_i}\xi(e_j)(\ch_j^*,\ch_i^*)=
-\la_j\xi(e_i)(e_{-1}^*,\ch_i^*)+\la_i\xi(e_j)(\ch_j^*,e_{-1}^*)\\
&=&-
\la_j\check{\al}_{-1,i}(e_i)-\la_i\check{\al}_{-1,j}(e_j),\\
0&=&\ad_{e_i}\xi(\ch_j)(\ch_i^*,e_j^*)-\ad_{\ch_j}\xi(e_i)(\ch_i^*,e_j^*)=
-\la_i\xi(\ch_j)(e_{-1}^*,e_j^*)+\la_j\xi(e_i)(e_{-1}^*,\ch_i^*)\\
&=&-\la_i\al_{-1,j}(\ch_j)+\la_j\check{\al}_{-1,i}(e_i),\\
0&=&\ad_{\ch_j}\xi(e_j)(\ch_j^*,e_j^*)-\ad_{e_j}\xi(\ch_j)(\ch_j^*,e_j^*)=
\la_j\xi(e_j)(\ch_j^*,e_{-1}^*)+\la_j\xi(\ch_j)(e_{-1}^*,e_j^*)\\
&=&-\la_j\check{\al}_{-1,j}(e_j)+\la_j\al_{-1,j}(\ch_j).
\end{eqnarray*}
We deduce that $\al_{-1,i}(\ch_i)=\check{\al}_{-1,i}(e_i)=0$ and hence
\begin{equation}\label{center}[e_{-1}^*,e_i^*]^*=\al_{-1,i}(e_{-1}) e_{-1}^*\quad\mbox{and}\quad [e_{-1}^*,\ch_i^*]^*=\check{\al}_{-1,i}(e_{-1}) e_{-1}^*.\end{equation}
To complete the proof, we will show that
$\al_{-1,i}(e_{-1})=\check{\al}_{-1,i}(e_{-1})=0$ and
$\al_{-1,0}=0$.

Since $\xi$ is a bialgebra structure, the dual bracket satisfies the Jacobi identity. Let us apply this identity to $e_{-1}^*,e_0^*,\ch_i^*$ for $i=1,\ldots,n$. We have
\begin{eqnarray*}
\;[[e_{-1}^*,e_0^*]^*,e_i^*]^*(e_i)&=&\al_{-1,0}(e_{-1})[e_{-1}^*,e_i^*]^*(e_i)+{\al}_{-1,0}(\ch_i)[\ch_i^*,
e_i^*]^*(e_i)\\&&+\sum_{j}{\al}_{-1,0}(e_j)[e_j^*, e_i^*]^*(e_i)+
\sum_{j\not=i}{\al}_{-1,0}(\ch_j)[\ch_j^*, e_i^*]^*(e_i)\\
&\stackrel{(\ref{center})}=&-{\al}_{-1,0}(\ch_i)c_{i,i}(e_i)+
\sum_{j}{\al}_{-1,0}(e_j)b_{j,i}(e_i)-\sum_{j\not=i}{\al}_{-1,0}(\ch_j)c_{i,j}(e_i)\\
&\stackrel{}=&-{\al}_{-1,0}(\ch_i)c_{i,i}(e_i)\quad\quad(\mbox{see Proposition \ref{pr2}}).\\
\;[[e_{i}^*,e_{-1}^*]^*,e_0^*]^*&\stackrel{(\ref{center})}=&-\al_{-1,i}(e_{-1})[e_{-1}^*,e_0^*]^*,\\
\;[[e_0^*,e_i^*]^*,e_{-1}^*]^*(e_i)&=&0.
\end{eqnarray*}The last equality is a consequence of the fact that $\xi(e_0)=0$ and (\ref{center}).
Hence, the Jacobi identity implies
$${\al}_{-1,0}(\ch_i)c_{i,i}(e_i)+\al_{-1,i}(e_{-1}){\al}_{-1,0}(e_i)=0.$$

Now
\begin{eqnarray*}
c_{i,i}(e_i)&=&\xi(e_i)(e_i^*, \ch_i^*)
\stackrel{(\ref{adjoint})}=-\xi(e_i)( \ad_{\ch_i}^*e_0^*,\ch_i^*)\\
&=&-\ad_{\ch_i}\xi(e_i)( e_0^*,\ch_i^*)
\stackrel{(\ref{cocycle})}=-\ad_{e_i}\xi(\ch_i)( e_0^*,\ch_i^*)\\
&\stackrel{(\ref{adjoint})}=&-\la_i\al_{-1,0}(\ch_i),\\
\al_{-1,i}(e_{-1})&=&\xi(e_{-1})(e_{-1}^*, e_i^*)
\stackrel{(\ref{adjoint})}=-\xi(e_{-1})(e_{-1}^*, \ad_{\ch_i}^*e_0^*)\\
&=&-\ad_{\ch_i}\xi(e_{-1})(e_{-1}^*, e_0^*)
\stackrel{(\ref{cocycle})}=-\left(\ad_{e_{-1}}\xi(\ch_i)(e_{-1}^*,
e_0^*)+\la_i\xi(e_i)(e_{-1}^*, e_0^*)\right)\\
&=&-\la_i\al_{-1,0}(e_i).\end{eqnarray*}

Hence
$({\al}_{-1,0}(\ch_i))^2+(\al_{-1,0}(e_i))^2=0,$ and then
${\al}_{-1,0}(\ch_i)=\al_{-1,0}(e_i)=0.$ Note that we have also shown that \begin{equation}\label{cii}c_{i,i}(e_i)=\al_{-1,i}(e_{-1})=0,\end{equation} and in the same way we can show that \begin{equation}\label{cii1}c_{i,i}(\ch_i)=\check{\al}_{-1,i}(e_{-1})=0.\end{equation}
To complete the proof, note that  $\al_{-1,0}(e_{-1})=0$ according to  Proposition \ref{pr4}.
\hfill $\square$
\begin{pr}\label{pr3}Let $\xi:\G_\la\too\G_\la\wedge\G_\la$ be a bialgebra
structure on $\G_\la$. Then, for $i,j=1,\ldots,n$ with $i\not=j$,
\begin{eqnarray*}
c_{i,i}&=&0,\\
b_{i,j}(e_i)&=&b_{i,j}(\ch_i)=b_{i,j}(e_j)=b_{i,j}(\ch_j)=0,\\
\check{b}_{i,j}(e_i)&=&\check{b}_{i,j}(\ch_i)=\check{b}_{i,j}(e_j)=\check{b}_{i,j}(\ch_j)=0,\\
c_{i,j}(e_i)&=&c_{i,j}(\ch_i)=c_{i,j}(e_j)=c_{i,j}(\ch_j)=0.\end{eqnarray*}

\end{pr}

{\bf Proof.} The vanishing of $c_{i,i}$ is a consequence of  Propositions \ref{pr1}-\ref{pr2} and (\ref{cii})-(\ref{cii1}).

Note that in Proposition \ref{pr2}, we have shown
\begin{eqnarray*}
b_{i,j}(e_i)&=&b_{i,j}(e_j)=
c_{i,j}(e_i)=c_{i,j}(\ch_j)=
\check{b}_{i,j}(\ch_i)=\check{b}_{i,j}(\ch_j)=0.\end{eqnarray*}
On the other hand, by using (\ref{cocycle}) and (\ref{adjoint}), we get
\begin{eqnarray*}
b_{i,j}(\ch_i)&=&\xi(\ch_i)(e_i^*,e_j^*)
=-\xi(\ch_i)(e_i^*,\ad_{\ch_j}^*e_0^*)\\
&=&-\ad_{\ch_j}\xi(\ch_i)(e_i^*,e_0^*)
=-\ad_{\ch_i}\xi(\ch_j)(e_i^*,e_0^*)\\
&=&-\la_i\xi(\ch_j)(e_{-1}^*,e_0^*)=0.\end{eqnarray*}
The same calculation gives the other relations.
\hfill$\square$\\

We end this section by proving Theorems \ref{main1}-\ref{flat}.\\

{\bf Proof of Theorem \ref{main1}.}\\

Suppose that $\xi(u)=\ad_u^\dag r+2e_0\wedge (J+\ad_{u_0})(u)$. Since $e_0$ is central and $J+\ad_{u_0}$ is a derivation, $\xi$ is a 1-cocycle with respect to the adjoint action. Remark that for any $u\in\G_\la$ and any $\al\in\G_\la$, $r(\ad_{u}^*e_0^*,\al)=i_{r_{\#}(\al)}\om(u).$ From this relation and by a straightforward computation one can show that the bracket on $\G_\la$ associated to $\xi$ is given by (\ref{bracketmain}). Moreover, $\xi$ defines a Lie bialgebra structure on $\G_\la$ iff this bracket satisfies Jacobi identity which is equivalent to
$$[e_0^*,[\al,\be]^*]^*+[\al,[\be,e_0^*]^*]^*+[\be,[e_0^*,\al]^*]^*=0,$$for all $\al,\be\in S^*$. Now $[e_0^*,[\al,\be]^*]^*=0$ and
\begin{eqnarray*}
 [\al,[\be,e_0^*]^*]^*&=&-2[\al, J^*\be]^*-[\al,i_{r_{\#}(\be)}\om]^*\quad\quad(J^*\be,i_{r_{\#}(\be)}\om\in S^*)\\
 &=&-2(\ad_{e_{-1}}^\dag r)(\al,J^*\be){e_{-1}^*}-(\ad_{e_{-1}}^\dag r)(\al,i_{r_{\#}(\be)}\om){e_{-1}^*}\\
 &=&-2(\ad_{e_{-1}}^\dag r)(\al,J^*\be){e_{-1}^*}-\om(r_{\#}(\be),(\ad_{e_{-1}}^\dag r)_{\#}(\al)){e_{-1}^*}.\end{eqnarray*}
  En conclusion the Jacobi identity holds iff
 $$\om((\ad_{e_{-1}}^\dag r)_{\#}(\al),r_{\#}(\be))+\om(r_{\#}(\al),(\ad_{e_{-1}}^\dag r)_{\#}(\be))
 -2(J^\dag\circ\ad_{e_{-1}}^\dag) r(\al,\be)=0$$which is equivalent to (\ref{boucetta}).

 Conversely, suppose that $\G_\la$ is generic and $\xi$ is a Lie bialgebra structure on $\G_\la$. By gathering all the relations shown in Propositions \ref{pr1}-\ref{pr3} and Lemma \ref{centrallemma} we get that  bracket on $\G_\la^*$ associated to $\xi$ is given by (the vanishing bracket are omitted)\begin{eqnarray*}
\;[e_0^*,e_i^*]^*&=&a_{0,i}(e_{-1})e_{-1}^*+\sum_{j=1}^n\left({a}_{0,i}(\ch_{j})\ch_j^*+
{a}_{0,i}(e_{j})e_j^*\right),\\
\;[e_0^*,\ch_i^*]^*&=&\check{a}_{0,i}(e_{-1})e_{-1}^*+\sum_{j=1}^n\left(\check{a}_{0,i}(\ch_{j})\ch_j^*+
\check{a}_{0,i}(e_{j})e_j^*\right),\\
\;[e_i^*,e_j^*]^*&=&\left(\la_j\check{a}_{0,j}(\ch_i)-\la_ia_{0,j}(e_i)\right)e_{-1}^*,\\
\;[e_i^*,\ch_j^*]^*&=&\left(-\la_i\check{a}_{0,j}(e_i)-\la_ja_{0,j}(\ch_i)\right)e_{-1}^*\\
&=&
\left(\la_ja_{0,i}(\ch_j)+\la_i\check{a}_{0,i}(e_j)\right)e_{-1}^*\quad i\not=j,\\
\;[\ch_i^*,\ch_j^*]^*&=&\left(\la_ja_{0,j}(e_i)-\la_i\check{a}_{0,j}(\ch_i)\right)e_{-1}^*.\end{eqnarray*}
The skew-symmetry of $[e_i^*,e_j^*]^*$ and $[\ch_i^*,\ch_j^*]^*$ and the two expressions of $[e_i^*,\ch_j^*]^*$ give
$$\left\{\begin{array}{ccc}
\la_j(\check{a}_{0,j}(\ch_i)-a_{0,i}(e_j))-\la_i(a_{0,j}(e_i)-\check{a}_{0,i}(\ch_j))&=&0,\\
-\la_i(\check{a}_{0,j}(\ch_i)-a_{0,i}(e_j))+\la_j(a_{0,j}(e_i)-\check{a}_{0,i}(\ch_j))&=&0,\\
\la_j(a_{0,i}(\ch_j)+a_{0,j}(\ch_i))+\la_i(\check{a}_{0,i}(e_j)+\check{a}_{0,j}(e_i))&=&0,\\
\la_i(a_{0,i}(\ch_j)+a_{0,j}(\ch_i))+\la_j(\check{a}_{0,i}(e_j)+\check{a}_{0,j}(e_i))&=&0.\end{array}\right.$$
Since the $\la_i$ are mutually distinct, these relations are equivalent to
\begin{equation}\label{symmetry}\check{a}_{0,j}(\ch_i)-a_{0,i}(e_j)=a_{0,j}(e_i)-\check{a}_{0,i}(\ch_j)=
a_{0,i}(\ch_j)+a_{0,j}(\ch_i)=\check{a}_{0,i}(e_j)+\check{a}_{0,j}(e_i)=0.\end{equation}
Put$$\left\{\begin{array}{lllclll}
r_0(e_0^*,e_i^*)&=&\frac1{\la_i}\check{a}_{0,i}(e_{-1}),&&r_0(e_0^*,\ch_i^*)&=&-\frac1{\la_i}{a}_{0,i}(e_{-1}),\\
r_0(e_i^*,e_j^*)&=&{a}_{0,i}(\ch_{j}),&&r_0(\ch_i^*,\ch_j^*)&=&-\check{a}_{0,i}(e_{j}),\\
r_0(e_i^*,\ch_j^*)&=&-{a}_{0,i}(e_{j}),&& i_{e_{-1}^*}r_0&=&0,\quad a_i=\check{a}_{0,i}(e_i),\end{array}\right.$$
and define $J$  by $Je_0=J_{e_{-1}}=0$, $Je_i=a_i\ch_i$ and $J\ch_i=-a_ie_i$. From (\ref{symmetry}),
$r_0\in\wedge^2\G_\la$. The endomorphism $J$ is a derivation which commutes with $\ad_{e_{-1}}$ and
   by comparing the brackets above to (\ref{exactbracket})  one can see that $\xi=\ad^\dag r_0+e_0\wedge J$. To complete the proof note that since $i_{e_{-1}^*}r_0=0$, there exists $r\in\wedge^2S$ and $u_0\in S$ such that $r_0=e_0\wedge u_0+r$.  \hfill $\square$\\

{\bf Proof of Theorem \ref{main2}.}
\begin{enumerate}\item Let $r\in\wedge^2\G_\la$ be a solution of the generalized classical Yang-Baxter equation. Then $\xi$ given by $\xi(u)=\ad_u^\dag r$ defines a Lie bialgebra structure on $\G_\la$. By Lemma \ref{centrallemma}, for any $i=1,\ldots,n$, $\ad_{e_{-1}}^\dag r(e_{-1}^*,e_i)=\ad_{e_{-1}}^\dag r(e_{-1}^*,\ch_i)=0$. These relations are equivalent to
$r(e_{-1}^*,e_i^*)=r(e_{-1}^*,\ch_i^*)=0$ and hence $$r=r(e_0^*,e_{-1}^*)e_0\wedge e_{-1}+e_0\wedge u_0+r_0,$$ where $r_0\in\wedge^2S$ and $u_0\in S$. Thus $$\xi(u)=\ad_u^\dag r=\ad_u^\dag r_0+2e_0\wedge(-\al\ad_{e_{-1}}-\frac12\ad_{u_0})(u),$$where $r(e_0^*,e_{-1}^*)=2\al$. Now according to Theorem \ref{main1}, $\xi$ defines a Lie bialgebra structure iff $r_0$ satisfies (\ref{boucetta}) with $J=-\al\ad_{e_{-1}}$ and this equivalent to $r_0$ satisfies (\ref{gyb1}).
\item Let $r\in\wedge^2\G_\la$ be a solution of the classical Yang-Baxter equation. In particular, $r$ is a solution of the generalized classical Yang-Baxter equation and from 1.,
    $$r=\al e_0\wedge e_{-1}+e_0\wedge u_0+r_0,$$ where $r_0\in\wedge^2S$ and $u_0\in S$. Now
    $$[r,r]=[r_0,r_0]+2\al e_0\wedge \ad^\dag_{e_{-1}}r_0+2e_0\wedge\ad_{u_0}^\dag r_0.$$Now, by using (\ref{bracket1}) and the formula (see \cite{dz})
\begin{equation}\label{zung}[r_0,r_0](\al,\be,\ga)=2\al\left([r_{0\#}(\be),r_{0\#}(\ga)]\right)+2\be\left([r_{0\#}(\ga),r_{0\#}(\al)]\right)+
2\ga\left([r_{0\#}(\al),r_{0\#}(\be)]\right)\end{equation}one can check easily that
$$[r_0,r_0]=2e_0\wedge\om_{r_0,r_0}.$$Thus $[r,r]=0$ iff
$$e_0\wedge\left(\om_{r_0,r_0}+\al\ad^\dag_{e_{-1}}r_0+\ad_{u_0}^\dag r_0\right)=0.$$
Since $r_0\in\wedge^2S$, $e_0\wedge \ad_{u_0}^\dag r_0=0$, $\om_{r_0,r_0},\ad^\dag_{e_{-1}}r_0\in\wedge^2S$ and hence the equation above is equivalent to (\ref{cyb1}).\hfill$\square$
\end{enumerate}

{\bf Proof of Theorem \ref{flat}}\\
\begin{enumerate}\item
We identify $(\G^*,[\;,\;]_r)$ to the Lie algebra  of left invariant vector fields on $G_r^*$  and $\prs^*$ to the restriction of $k^*$ to $\G^*$. From the definition $[\;,\;]_r$ given by (\ref{dualbracket1}) and since $\ad_{u}^*$ is skew-symmetric for any $u\in\G$, one can see easily that, for any $\al,\be\in\G^*$,
$$\na^*_\al\be=-\ad_{r_\#(\al)}^*\be.$$Thus the curvature of $\na^*$ is given by
$$R(\al,\be)\ga=\na^*_{[\al,\be]_r}\ga-[\na^*_\al,\na^*_\be]\ga=\ad^*_{[r_\#(\al),r_\#(\be)]-r_\#([\al,\be]_r)}\ga.$$
Now, from  formula (\ref{zung}) one can  deduce that
$$\ga([r_\#(\al),r_\#(\be)]-r_\#([\al,\be]_r))=\frac12[r,r](\al,\be,\ga).$$
From this relation, we deduce that if $r$ is a solution of (\ref{cyb}) then $R$ vanishes identically. We return to the general case and
we denote by $u_r(\al,\be)$ the element of $\G$ defined by
\begin{equation}\label{ls}\ga(u_r(\al,\be))=\frac12[r,r](\al,\be,\ga).\end{equation}
Thus
$$R(\al,\be)\ga=\ad^*_{u_r(\al,\be)}\ga.$$Let us compute $\na^* R$.
We have, for any $\al,\be,\ga,\rho\in\G^*$,
\begin{eqnarray*}
(\na_{\rho}^*  R)(\al,\be,\ga)&=&\na_{\rho}^*  (R(\al,\be,\ga))- R(\na_{\rho}^*\al,\be,\ga)-
 R(\al,\na_{\rho}^*\be,\ga)- R(\al,\be,\na_{\rho}^*\ga)\\
&=&-\ad_{r_\#(\rho)}^*\circ\ad^*_{u_r(\al,\be)}\ga+
\ad^*_{u_r(\ad_{r_\#(\rho)}^*\al,\be)}\ga+\ad^*_{u_r(\al,\ad_{r_\#(\rho)}^*\be)}\ga\\&&+
\ad^*_{u_r(\al,\be)}\circ\ad_{r_\#(\rho)}^*\ga\\
&=&\ad^*_{[r_\#(\rho),u_r(\al,\be)]}\ga+
\ad^*_{u_r(\ad_{r_\#(\rho)}^*\al,\be)}\ga+\ad^*_{u_r(\al,\ad_{r_\#(\rho)}^*\be)}\ga\\
&=&\ad^*_{w(\al,\be)}\ga,
\end{eqnarray*}where
$$w(\al,\be)=[r_\#(\rho),u_r(\al,\be)]+u_r(\ad_{r_\#(\rho)}^*\al,\be)+u_r(\al,\ad_{r_\#(\rho)}^*\be).$$
Now, for any $\mu\in\G^*$, we have
\begin{eqnarray*}
\mu\left(w(\al,\be)\right)&=&\ad_{r_\#(\rho)}^*\mu(u_r(\al,\be))+\mu(u_r(\ad_{r_\#(\rho)}^*\al,\be))
+\mu(u_r(\al,\ad_{r_\#(\rho)}^*\be))\\
&\stackrel{(\ref{ls})}=&\frac12[r,r](\al,\be,\ad_{r_\#(\rho)}^*\mu)+
\frac12[r,r](\ad_{r_\#(\rho)}^*\al,\be,\mu)+\frac12[r,r](\al,\ad_{r_\#(\rho)}^*\be,\mu)\\
&=&\frac12(\ad_{r_\#(\rho)}[r,r])(\al,\be,\mu)\stackrel{(\ref{gyb})}=0.\end{eqnarray*}
This achieves to show that $\na^*R=0$ and the first part of the theorem follows.

  \item According to a result of Aubert and Medina (see \cite{medina}), a flat left invariant pseudo-Riemannian metric on a Lie group is complete if and only if this group is unimodular and in this case the group is solvable.  \hfill $\square$\end{enumerate}

\section{Solutions of (\ref{boucetta}), (\ref{gyb1}) and  (\ref{cyb1}) when $\dim G_\la\leq6$}\label{example}
In this section, We give all the bialgebras structures, all the solutions of generalized classical Yang-Baxter equation and all the solutions of the classical Yang-Baxter equation on $\G_\la$ when $\G_\la$ is generic and $\dim\G_\la\leq 6$. When $\dim\G_\la\geq 8$, we give a large class of solutions.
According to Theorems \ref{main1}-\ref{main2}, we will solve   (\ref{boucetta}), (\ref{gyb1}) and  (\ref{cyb1}).\\
Note first that if   $J:\G_\la\too\G_\la$ is derivation  commuting with $\ad_{e_{-1}}$ and satisfying $J(e_{-1})=J(e_0)=0$ then there exists $a=(a_1,\ldots,a_n)\in\R^n$ such that for any $i=1,\ldots,n$ such that
$$J(e_i)=a_i\ch_i\quad\mbox{and}\quad J(\ch_i)=-a_ie_i.$$We shall denote by $J_a$ such an endomorphism.\\
When $\dim\G_\la=4$, the situation is simple.  We have $\B=\{e_{-1},e_0,e_1,\ch_1\}$, $S=span\{e_1,\ch_1\}$ and
$\wedge^2S=span\{e_1\wedge\ch_1\}$. Moreover, $\ad_{e_{-1}}^\dag(e_1\wedge\ch_1)=J_a^\dag(e_1\wedge\ch_1)=0$ and
$\om_{e_1\wedge\ch_1,e_1\wedge\ch_1}=e_1\wedge \ch_1$. By applying Theorems \ref{main1}-\ref{main2}, we get the following proposition.

 \begin{pr}\label{dim4} Let $\la\in\R$ and let $\G_\la$ be the associated 4-dimensional oscillator Lie algebra. Then:
 \begin{enumerate}\item $\xi:\G_\la\too\wedge^2\G_\la$ defines a Lie bialgebra structure on $\G_\la$ if and only if there exists $a,\al\in\R$ and $u_0\in S$ such that
 $$\xi(u)=\al \ad_u^\dag(e_1\wedge\ch_1)+e_0\wedge(J_a+\ad_{u_0})(u).$$
 \item A bivector $r\in\wedge^2\G_\la$ is a solution of the generalized classical Yang-Baxter equation if and only if $r=e_0\wedge u+\al e_1\wedge\ch_1,$ where $\al\in\R$ and $u\in\G_\la$.
     \item A bivector $r\in\wedge^2\G_\la$ is a solution of the  classical Yang-Baxter equation if and only if $r=e_0\wedge u,$ where $\al\in\R$ and $u\in\G_\la$.

 \end{enumerate}\end{pr}

 We return now to the general case. Consider a generic $2n+1$-dimensional oscillator Lie algebra $\G_\la$ and fix $a=(a_1,\ldots,a_n)\in\R^n$. For any $1\leq i,j\leq n$, put
\begin{eqnarray*}
r_{ij}&=&e_i\wedge e_j,\;
\check{r}_{ij}=\ch_i\wedge \ch_j,\;
s_{ij}=e_i\wedge\ch_j,\;
\check{s}_{ij}=\ch_i\wedge e_j,\;
t_i=e_i\wedge \ch_i.\end{eqnarray*}
A direct computation gives\begin{eqnarray*}
\om_{r_{ij},r_{ij}}&=&\om_{\check{r}_{ij},\check{r}_{ij}}=\om_{s_{ij},s_{ij}}=
\om_{\check{s}_{ij},\check{s}_{ij}}=0,\\
\om_{r_{ij},s_{ij}}&=&\om_{r_{ij},\check{s}_{ij}}=\om_{\check{r}_{ij},s_{ij}}=\om_{\check{r}_{ij},\check{s}_{ij}}=0,\\
\om_{r_{ij},\check{r}_{ij}}&=&=-\om_{s_{ij},\check{s}_{ij}}=
\frac12(t_i+t_j),\\
\om_{t_i,r_{ij}}&=&\frac12r_{ij},\;\om_{t_i,\check{r}_{ij}}=\frac12\check{r}_{ij},\;
\om_{t_i,s_{ij}}=\frac12s_{ij},\;\om_{t_i,\check{s}_{ij}}=\frac12\check{s}_{ij}\\
\om_{t_j,r_{ij}}&=&\frac12r_{ij},\;\om_{t_j,\check{r}_{ij}}=\frac12\check{r}_{ij},\;
\om_{t_j,s_{ij}}=\frac12s_{ij},\;\om_{t_j,\check{s}_{ij}}=\frac12\check{s}_{ij},\\
\om_{t_i,t_i}&=&t_i,\;\om_{t_i,t_j}=0.\end{eqnarray*}

Fix $1\leq i<j\leq n$ and let us find the solutions of (\ref{boucetta}), (\ref{gyb1}) and (\ref{cyb1}) of the form $r=p_{ij}+c_it_i+c_jt_j$ where $p_{ij}\in span\{r_{ij},\check{r}_{ij},s_{ij},\check{s}_{ij}\}$. In order to simplify the computations, let us introduce the following new basis:
$$E_{ij}=s_{ij}+\check{s}_{ij},\;\check{E}_{ij}=-r_{ij}+\check{r}_{ij},\;F_{ij}=-s_{ij}+\check{s}_{ij},\;
\check{F}_{ij}=r_{ij}+\check{r}_{ij}.$$
We have
$$J_a^\dag(E_{ij})=(a_i+a_j)\check{E}_{ij},\;J_a^\dag(\check{E}_{ij})=-(a_i+a_j){E}_{ij},\;
J_a^\dag(F_{ij})=(a_j-a_i)\check{F}_{ij},\;J_a^\dag(\check{F}_{ij})=(a_i-a_j){F}_{ij}.$$
Since $\ad_{e_{-1}}^\dag=J_\la^\dag$, similar relations holds for $\ad_{e_{-1}}^\dag$. On the other hand, one can see easily that $(E_{ij},\check{E}_{ij},F_{ij},\check{F}_{ij})$ is $\om$-orthogonal and
$$\om_{E_{ij},E_{ij}}=\om_{\check{E}_{ij},\check{E}_{ij}}=-\om_{F_{ij},F_{ij}}=
-\om_{\check{F}_{ij},\check{F}_{ij}}=-(t_i+t_j).$$
Let $(a,\check{a},b,\check{b})$ be the coordinates of $p_{ij}$ in $(E_{ij},\check{E}_{ij},F_{ij},\check{F}_{ij})$. From the relations above we get\begin{eqnarray*}
\om_{r,r}+\al \ad_{e_{-1}}^\dag r&=&c_i^2t_1+c_j^2t_j+(c_i+c_j)p_{ij}+(b^2+\check{b}^2-a^2-\check{a}^2)(t_i+t_j)\\
&&+\al(\la_i+\la_j)\left(-\check{a}E_{ij}+a\check{E}_{ij}\right)\\
&&+\al(\la_j-\la_i)\left(-\check{b}F_{ij}
+b\check{F}_{ij}\right),\\
\om_{r,ad_{e_{-1}}^\dag r}-(J_a^\dag\circ \ad_{e_{-1}}^\dag)r&=&\frac12(c_i+c_j)(\la_i+\la_j)\left(-\check{a}E_{ij}+a\check{E}_{ij}\right)\\
&&+\frac12(c_i+c_j)(\la_j-\la_i)\left(-\check{b}F_{ij}
+b\check{F}_{ij}\right)\\
&&+(\la_i+\la_j)(a_i+a_j)\left(a{E}_{ij}+\check{a}\check{E}_{ij}\right)
\\&&+(\la_j-\la_i)(a_j-a_i)\left(b{F}_{ij}+\check{b}\check{F}_{ij}\right),\\
\om_{r,\ad_{e_{-1}}^\dag r}+\al (\ad_{e_{-1}}^\dag\circ \ad_{e_{-1}}^\dag)r&=&
\frac12(c_i+c_j)(\la_i+\la_j)\left(-\check{a}E_{ij}+a\check{E}_{ij}\right)\\
&&+\frac12(c_i+c_j)(\la_j-\la_i)\left(-\check{b}F_{ij}
+b\check{F}_{ij}\right)\\
&&-\al(\la_i+\la_j)^2\left(a{E}_{ij}+\check{a}\check{E}_{ij}\right)
\\&&-\al(\la_j-\la_i)^2\left(b{F}_{ij}+\check{b}\check{F}_{ij}\right).
\end{eqnarray*}

\begin{enumerate}\item\begin{enumerate}\item For $\al\not=0$, $r=p_{ij}+c_it_i+c_jt_j$ is a solution of (\ref{cyb1}) if and only if $r=0$.\item
For $\al=0$, $r=p_{ij}+c_it_i+c_jt_j$ is a solution of (\ref{cyb1}) if and only if $r=aE_{ij}+\check{a}\check{E}_{ij}+bF_{ij}+\check{b}\check{F}_{ij}+c(t_i-t_j)$ and
$c^2=a^2+\check{a}^2-b^2-\check{b}^2$.\end{enumerate}
\item\begin{enumerate}\item For $\al=0$, $r=p_{ij}+c_it_i+c_jt_j$ is a solution of (\ref{gyb1}) if and only if $r=c_it_i+c_jt_j$ or $r=aE_{ij}+\check{a}\check{E}_{ij}+bF_{ij}+\check{b}\check{F}_{ij}+c(t_i-t_j)$.
   \item  For $\al\not=0$, $r=p_{ij}+c_it_i+c_jt_j$ is a solution of (\ref{gyb1}) if and only if $r=c_it_i+c_jt_j$.
   \end{enumerate}
   \item \begin{enumerate}\item If $a_1+a_2\not=0$ and $a_2-a_1\not=0$ then $r=p_{ij}+c_it_i+c_jt_j$ is a solution of (\ref{boucetta}) if and only if $r=c_it_i+c_jt_j$.
   \item If $a_1+a_2=0$ and $a_2-a_1\not=0$ then $r=p_{ij}+c_it_i+c_jt_j$ is a solution of (\ref{boucetta}) if and only if $r=c_it_i+c_jt_j$ or $r=aE_{ij}+\check{a}\check{E}_{ij}+c(t_i-t_j)$.
   \item If $a_1+a_2\not=0$ and $a_2-a_1=0$ then $r=p_{ij}+c_it_i+c_jt_j$ is a solution of (\ref{boucetta}) if and only if $r=c_it_i+c_jt_j$ or $r=bF_{ij}+\check{b}\check{F}_{ij}+c(t_i-t_j)$.
       \item If $a_1=a_2=0$ then $r=p_{ij}+c_it_i+c_jt_j$ is a solution of (\ref{boucetta}) if and only if $r=c_it_i+c_jt_j$ or $r=aE_{ij}+\check{a}\check{E}_{ij}+bF_{ij}+\check{b}\check{F}_{ij}+c(t_i-t_j)$.

   \end{enumerate}

\end{enumerate}

Remark that if $r_1=p_{ij}+c_it_i+c_jt_j$ and $r_2=p_{lk}+c_lt_l+c_kt_k$ are two solutions of (\ref{boucetta}), (\ref{gyb1}) or (\ref{cyb1}) and $\{i,j\}\cap\{l,k\}=\emptyset$ then $r_1+ r_2$ is a solution of (\ref{boucetta}), (\ref{gyb1}) or (\ref{cyb1}). So we can construct a large class of solutions. When $\dim\G_\la=6$ we have constructed all the solutions. Let us summarize this case.
\begin{pr}\label{dim6} Let $\la=(\la_1,\la_2)$ with $\la_1<\la_2$ and let $\G_\la$ be the associated 6-dimensional oscillator Lie algebra. Then:
 \begin{enumerate}\item $\xi:\G_\la\too\wedge^2\G_\la$ defines a Lie bialgebra structure on $\G_\la$ if and only if $\xi$ has one of the following forms:
 \begin{enumerate}\item
 $\xi(u)= \ad_u(c_1t_1+c_2t_2)+e_0\wedge(J_a+\ad_{u_0})(u),$ where $c_1,c_2\in\R$,  $a\in\R^2$ and $u_0\in S$;
 \item $\xi(u)= \ad_u(dE_{12}+\check{d}\check{E}_{12}+c(t_1-t_2))+e_0\wedge(J_{(a,-a)}+\ad_{u_0})(u),$ where $a,d,\check{d},c\in\R$ and $u_0\in S$;
 \item $\xi(u)= \ad_u(dF_{12}+\check{d}\check{F}_{12}+c(t_1-t_2))+e_0\wedge(J_{(a,a)}+\ad_{u_0})(u),$ where $a,d,\check{d},c\in\R$ and $u_0\in S$;
     \item $\xi(u)= \ad_u(dE_{12}+\check{d}\check{E}_{12}+bF_{12}+\check{b}\check{F}_{12}+c(t_1-t_2))+e_0\wedge(\ad_{u_0})(u),$ where $b,\check{b},d,\check{d},c\in\R$ and $u_0\in S$;
 \end{enumerate}
 \item A bivector $r\in\wedge^2\G_\la$ is a solution of the generalized classical Yang-Baxter equation if and only if  $r=e_0\wedge u+c_1t_1+c_2t_2$, where $u\in\G_\la$ and $c_1,c_2\in\R$,
     or $r=e_0\wedge u+aE_{12}+\check{a}\check{E}_{12}+bF_{12}+\check{b}\check{F}_{12}+c(t_1-t_2)$, where $u\in S$ and $a,\check{a},b,\check{b},c\in\R$.
     \item A bivector $r\in\wedge^2\G_\la$ is a solution of the  classical Yang-Baxter equation if and only if $r=e_0\wedge u,$ where $\al\in\R$ and $u\in\G_\la$, or $r=e_0\wedge u+
         aE_{12}+\check{a}\check{E}_{12}+bF_{12}+\check{b}\check{F}_{12}+c(t_1-t_2)$ where $u\in S$ and
         $c^2=a^2+\check{a}^2-b^2-\check{b}^2$.

 \end{enumerate}\end{pr}
 By using Corollary \ref{co} and the proposition above, one can prove easily the following result.
 \begin{pr} Let $\la=(\la_1,\la_2)$ with $\la_1<\la_2$ and let $\G_\la$ be the associated 6-dimensional oscillator Lie algebra. Then for any solution $r$ of the classical Yang-Baxter equation  $(\G_\la^*,[\;,\;]_r)$ is unimodular.\end{pr}

\section{Examples}\label{example1}
\begin{enumerate}\item
 By using Theorem \ref{flat}, we will  build an example of  6-dimensional Lie groups endowed with a complete left invariant flat Lorentzian metric.\\
According to Proposition \ref{dim6},
$$r=e_0\wedge e_1+ e_1\wedge \ch_2+\ch_1\wedge e_2+e_1\wedge\ch_1-e_2\wedge\ch_2$$ is a solution of (\ref{cyb}) on $\G_\la$ ($\la=(\la_1,\la_2)$). By using (\ref{exactbracket}), we can check that the bracket $[\;,\;]_r$ on $\G_\la^*$ associated to $r$ is given by
\begin{eqnarray*}
\;[e_0^*,e_1^*]_r&=&-e_1^*-e_2^*,
\;[e_0^*,e_2^*]_r=e_1^*+e_2^*,
\;[e_0^*,\ch_1^*]_r=\la_1e_{-1}^*-\ch_1^*+\ch_2^*,\\
\;[e_0^*,\ch_2^*]_r&=&-\ch_1+\ch_2,\;[e_2^*,\ch_2^*]_r=[e_1^*,\ch_1^*]_r=[e_1^*,\ch_2^*]_r=[e_2^*,\ch_1^*]_r=0,\\
\;[e_1^*,e_2^*]_r&=&-[\ch_1^*,\ch_2^*]_r=-(\la_1+\la_2)e_{-1}^*.
\end{eqnarray*}
Remark that ${\cal I}=span\{e_{-1}^*,e_1^*,e_2^*,\ch_1^*,\ch_2^*\}$ is an ideal of $(\G_\la^*,[\;,\;]_r)$ isomorphic to
the 5-dimensional Heisenberg Lie algebra and
$(\G_\la^*,[\;,\;]_r)$ is isomorphic to the semi-direct product of $\R e_0^*$ with ${\cal I}$ and the action is given by the restriction of $\ad_{e^*_0}$ to  $\cal I$. Let $\textbf{k}_\la$ be the adjoint invariant Lorentzian product on $\G_\la$ given by (\ref{kl}). The corresponding  symmetric bilinear form $\textbf{k}_\la^*$  is entirely determined by the relations
$$\textbf{k}_\la^*(e_0^*,e_{-1}^*)=1,\;\textbf{k}_\la^*(e_i^*,e_i^*)=\textbf{k}_\la^*(\ch_i^*,\ch_i^*)=\la_i,\; i=1,2$$and induces, according to Theorem \ref{flat},
a complete flat  left invariant Lorentzian metric on the connected and simply connected Lie group $G_\la^*$.

\item Let $(\G,k)$ be an orthogonal Lie algebra an $r:\G^*\too\G$ a solution of (\ref{cyb}). If $\mbox{Im}r$ is $k$-nondegenerate then $\G_r^*$ is a solvable Lie algebra. Indeed, under the hypothesis, $\mbox{Im}r$ is a nilpotent Lie algebra (see \cite{mr}) and then $\G_r^*$ is an extension of a nilpotent  algebra by an Abelian
algebra.
Nevertheless, the following example  shows that $\G_r^*$ can be solvable with $\mbox{Im}r$ degenerate.

Let $\G=sl(2,\R)$ and let $\B=\{e_1,e_2,e_3\}$ the basis of $\G$ where $$e_1=\left(\begin{array}{cc}1&0\\0&-1\end{array}\right),\;\;
e_2=\left(\begin{array}{cc}0&1\\0&0\end{array}\right)\quad\mbox{ and} \quad e_3=\left(\begin{array}{cc}0&0\\-1&0\end{array}\right).$$ We have
$$[e_1,e_2]=2e_2,\;\;[e_1,e_3]=-2e_3\quad\mbox{and}\quad[e_2,e_3]=-e_1.$$

The symmetric 2-form
$$k(a,b)=\tr{(ab)}$$ is an orthogonal structure on $\G$. The matrix of $k$ and $k^*$ in $\B$ and $\B^*$ are given by
$$\mathrm{M}(k,\B)=\left(\begin{array}{ccc}2&0&0\\0&0&-1\\0&-1&0\end{array}\right)\quad
\mbox{and}\quad\mathrm{M}(k^*,\B^*)=\left(\begin{array}{ccc}\frac12&0&0\\0&0&-1\\0&-1&0\end{array}\right).$$
Let $r:\G^*\too\G$ be a linear endomorphism which is skew-symmetric, i.e., $\al(r(\be))=-\be(r(\al))$ for any $\al,\be\in\G^*$. Denote by $\left(\begin{array}{ccc}0&a&b\\-a&0&c\\-b&-c&0\end{array}\right)$ the matrix of $r$ in the basis $\B^*$ and $\B$. We have
\begin{eqnarray*}
\;[r(e_1^*),r(e_2^*)]&=&[-ae_2-be_3,ae_1-ce_3]\\
&=&a^2[e_1,e_2]+ac[e_2,e_3]+ab[e_1,e_3]\\
&=&-ace_1+2a^2e_2-2abe_3,\\
\;[r(e_2^*),r(e_3^*)]&=&[ae_1-ce_3,be_1+ce_2]\\
&=&ac[e_1,e_2]+bc[e_1,e_3]+c^2[e_2,e_3]\\
&=&-c^2e_1+2ace_2-2bce_3,\\
\;[r(e_3^*),r(e_1^*)]&=&[be_1+ce_2,-ae_2-be_3]\\
&=&-ab[e_1,e_2]-b^2[e_1,e_3]-bc[e_2,e_3]\\
&=&bce_1-2abe_2+2b^2e_3.
\end{eqnarray*}The endomorphism $r$ is a solution of (\ref{cyb}) if and only if
$$e_3^*\left([r(e_1^*),r(e_2^*)]\right)+e_1^*\left([r(e_2^*),r(e_3^*)]\right)
+e_2^*\left([r(e_3^*),r(e_1^*)]\right)=0.$$This equation is clearly equivalent to
\begin{equation}\label{eq1}4ab+c^2=0.\end{equation}
Let $r$ be a solution of (\ref{cyb}). It induces on $\G^*$ a Lie bracket given by (\ref{dualbracket1}).
A direct computation gives
$$[e_1^*,e_2^*]_r=-2ae_1^*-ce_2^*,\;[e_1^*,e_3^*]_r=2be_1^*-ce_3^*,\;[e_2^*,e_3^*]_r=2be_2^*+2ae_3^*.$$
One can see easily that if $r$ is non zero then  the derived ideal $[\G^*,\G^*]_r$ is abelian and of dimension 2 and hence $\G^*$ is isomorphic to a semi-direct product of a plan by a line where the action of the line is given by the identity.\\
According to Theorem \ref{flat}, the connected and simply connected Lorentzian Lie group associated to $(\G^*,[\;,\;]_r,k^*)$ is flat and non complete.
\end{enumerate}

{\bf Acknowledgments}\\

This research was conducted within the framework of Action concert\'ee CNRST-CNRS.\\
 This paper was partially written during a stay of the second author as a guest at the University of Antioquia in Colombia.

\noindent Mohamed Boucetta\\
Facult\'e des Sciences et Techniques \\
BP 549 Marrakech, Morocco.
\\
Email: {\it boucetta@fstg-marrakech.ac.ma }
\bigskip

\noindent Alberto Medina\\ Universit\'e Montpellier 2 \\Case Courrier 051,
UMR CNRS 5149\\
Place Eug\`ene Bataillon 34095\\ MONTPELLIER Cedex France\\
Email: {\it medina@math.univ-montp2.fr}

\end{document}